\documentclass[leqno]{article}
\usepackage{amssymb, amsmath, amsfonts}

\newtheorem{theorem}{Theorem}[section]

\newtheorem{lemma}[theorem]{Lemma}
\newtheorem{proposition}[theorem]{Proposition}
\newtheorem{definition}[theorem]{Definition}
\newtheorem{remark}[theorem]{Remark}

\numberwithin{equation}{section}

\def\sqr#1#2{{\vcenter{\vbox{\hrule height.#2pt
    \hbox{\vrule width.#2pt height#1pt \kern#1pt
    \vrule width.#2pt}
    \hrule height.#2pt}}}}

\def\ga{\gamma}

\def\De{\Delta}
\def\de{\delta}

\def\la{\lambda}
\def\eps{\varepsilon}

\def\vka{\varkappa}

\def\om{\omega}
\def\Om{\Omega}

\def\pa{\partial}

\def\diam{\hbox{diam}}

\def\bB{\bar{B}}

\def\tB{\tilde B}

\def\tQ{\tilde Q}

\font\bbb=msbm10

\def\R{\hbox{{\bbb R}}}

\def\cL{{\cal L}}

\def\capa{\, {\rm cap}\, }
\def\mes{\,{\rm mes}\, }
\def\Lip{\, {\rm Lip}}

\def\cN{{\cal N}}

\def\cG{{\mathcal G}}
\def\bcG{{\bar{\mathcal G}}}
\def\cE{{\mathcal E}}
\def\cD{{\mathcal D}}

\def\cV{\mathbb V}

\def\cbcG{\complement\bcG}

\def\msni{\medskip\noindent}

\def\ms{\medskip}
\def\bs{\bigskip}

\newdimen\theight
\def\TeXref#1{%
          \leavevmode\vadjust{\setbox0=\hbox{{\tt
                  \quad\quad  {\small \textrm #1}}}%
          \theight=\ht0
          \advance\theight by \lineskip
          \kern -\theight \vbox to
          \theight{\rightline{\rlap{\box0}}%
          \vss}%
          }}%

\begin{document}
\title{Discreteness of spectrum and positivity criteria 
for Schr\"odinger operators}
\author{
{\bf Vladimir Maz'ya}
\thanks{Research partially supported by the Department of
Mathematics and the Robert G. Stone Fund at Northeastern
University}
\\
Department of Mathematics\\
The Ohio State  University\\
Columbus, OH 43210,
USA\medskip\\
E-mail: vlmaz@mai.liu.se
\bigskip\\
{\bf Mikhail Shubin
\thanks{Research partially supported by NSF grant DMS-0107796} }
\\
Department of Mathematics\\
Northeastern University\\
Boston, MA 02115,
USA
\medskip\\
E-mail: shubin@neu.edu
}

\date{}

\maketitle

\begin{abstract}

We provide a  class of necessary and sufficient conditions
for the discreteness of spectrum of  Schr\"odinger operators with 
scalar potentials which are semibounded below.  
The classical discreteness of spectrum criterion by 
A.M.Molchanov (1953) uses a notion of negligible set
in a cube as a set whose Wiener's capacity is  
less than a small constant times the capacity of the cube.
We prove that this constant can be taken arbitrarily between
0 and 1. This solves a problem formulated  by I.M.Gelfand in 1953.
Moreover, we  extend the notion of negligibility by allowing 
the constant to depend on the size of the cube. We give a complete description
of all negligibility conditions of this kind.
The a priori equivalence of our conditions involving 
different negligibility classes 
is a non-trivial property of the  capacity.
We also establish similar strict positivity criteria for the Schr\"odinger operators
with non-negative potentials.

\end{abstract}

\section{Introduction}\label{S:Intro}

In 1934, K.~Friedrichs \cite{Friedrichs} proved that the spectrum of the Schr\"odinger operator 
$-\De+V$ in $L^2(\R^n)$ with a locally integrable potential $V$ is discrete provided
$V(x)\to+\infty$ as $|x|\to\infty$ (see also \cite{Berezin-Shubin, Reed-Simon}). 
On the other hand,  if we assume that $V$ is semi-bounded below, then
the discreteness of spectrum
easily  implies that for every
$d>0$
\begin{equation}\label{E:ness-simple}
\int_{Q_d}V(x)dx\to +\infty \quad {\rm as} \quad Q_d\to\infty,
\end{equation}
where $Q_d$ is an open cube with the edge length $d$ and with the edges parallel to coordinate axes,
$Q_d\to\infty$ means that the cube $Q_d$ goes to infinity (with fixed $d$). This was first noticed by
A.M.Molchanov in 1953 (see \cite{Molchanov}) who also showed that this condition is in fact
necessary and sufficient in case $n=1$ but not sufficient for $n\ge 2$. 
Moreover, in the same paper 
Molchanov discovered a modification of condition \eqref{E:ness-simple} which is fully
equivalent to the discreteness of spectrum in the case  $n\ge 2$. 
It states that for every $d>0$
\begin{equation}\label{E:main-cond}
\inf_F\int_{Q_d\setminus F} V(x) dx \to +\infty \quad {\rm as} \quad Q_d\to\infty,
\end{equation}
where infimum is taken over all compact subsets $F$ of the closure $\bar Q_d$ which are called
{\it negligible}.  The negligibility of $F$ in the sense of Molchanov means that
$\capa(F)\le\ga\capa(Q_d)$, where $\capa$ is the Wiener capacity and $\ga>0$ 
is a sufficiently small constant. More precisely, Molchanov proved that we can take $\ga=c_n$ 
where for $n\ge 3$ 
\begin{equation*}\label{E:cn}
c_n=(4n)^{-4n}(\capa(Q_1))^{-1}.
\end{equation*}
Proofs of Molchanov's result can be found also in \cite{Mazya8, Edmunds-Evans, Kondrat'ev-Shubin}.
In particular, the books \cite{Mazya8, Edmunds-Evans} contain a proof which first appeared in
\cite{Mazya73} and is different from the original Molchanov proof. We will not list numerous
papers related to the discreteness of spectrum conditions for one- and multidimensional 
Schr\"odinger operators. Some references can be found in \cite{Mazya8, Kondrat'ev-Shubin, KMS}.

As early as in 1953, I.M.Gelfand raised the question about the best possible constant $c_n$
(personal communication). 
In this paper we answer this question by proving 
that $c_n$ can be replaced by an arbitrary constant $\gamma$, $0<\gamma<1$.   

We even establish a stronger result.  We allow negligibility conditions of the form
\begin{equation}\label{E:negligibility}
\capa(F)\le\ga(d)\capa(Q_d)
\end{equation}
and completely describe all admissible functions $\gamma$. More precisely,
in the necessary condition for the discreteness of spectrum we allow arbitrary functions 
$\ga:(0,+\infty)\to (0,1)$. In the sufficient condition we can admit
arbitrary functions $\gamma$ with values in $(0,1)$, defined for $d>0$ in a neighborhood 
of $d=0$ and satisfying
\begin{equation}\label{E:gamma-cond1}
\limsup_{d\downarrow 0} d^{-2}\ga(d) = +\infty.
\end{equation}
On the other hand, if $\ga(d)=O(d^2)$ 
in the negligibility condition \eqref{E:negligibility}, then the condition 
\eqref{E:main-cond} is no longer sufficient, i.e. it may happen that 
it is satisfied but the spectrum  is not  discrete. 

All conditions \eqref{E:main-cond} involving functions $\ga:(0,+\infty)\to (0,1)$,
satisfying \eqref{E:gamma-cond1}, are necessary and sufficient for the discreteness of spectrum. 
Therefore two conditions with different functions $\gamma$ are equivalent, which is far from being 
obvious a priori. This equivalence means the following striking effect: if \eqref{E:main-cond}
holds for very small sets $F$, then it also holds for sets $F$ which almost fill the corresponding cubes.

Another important question is whether the operator $-\De+V$ with $V\ge 0$ is strictly positive, i.e.
the spectrum is separated from $0$. Unlike the discreteness of spectrum conditions, 
it is the  large values of $d$ which are relevant here. The following necessary
and sufficient condition  for the strict positivity
was obtained in \cite{Mazya73} (see also \cite{Mazya8}, Sect.12.5): 
there exist positive constants
$d$ and $\varkappa$ such that for all cubes $Q_d$
\begin{equation}\label{E:pos-cond}
\inf_F\int_{Q_d\setminus F}V(x)dx \ge \varkappa\;,
\end{equation}
where the infimum is taken over all compact sets $F\subset \bar Q_d$ which are negligible 
in the sense of Molchanov. We prove that here again 
an arbitrary constant $\ga\in (0,1)$ in the negligibility condition 
\eqref{E:negligibility} is admissible. 

The above mentioned results are proved in this paper in a more
general context.  
The family of cubes $Q_d$ is
replaced  by a family of arbitrary bodies homothetic to a standard bounded domain which is
star-shaped with respect to a ball.
Instead of locally integrable potentials $V\ge 0$ we consider
positive measures. 
We also include operators in arbitrary open subsets of $\R^n$ 
with the Dirichlet boundary conditions.

\section{Main results}\label{S:main}

Let $\cV$ be a positive Radon measure in  an open set $\Omega\subset\R^n$. We will consider
the Schr\"odinger operator which is formally given by an expression $-\De+\cV$. It
is defined in $L^2(\Om)$ by the quadratic form
\begin{equation}\label{E:form}
h_{\cV}(u,u)=\int_\Om|\nabla u|^2dx+\int_{\Om}|u|^2\cV(dx), \quad u\in C_0^\infty(\Om),
\end{equation}
where $C_0^\infty(\Om)$ is the space of all $C^\infty$-functions 
with compact support in $\Om$. For the associated operator to be well defined 
we need a closed form.  
The form above is closable in $L^2(\Om)$  if and only if $\cV$ is absolutely continuous 
with respect to the Wiener capacity, i.e. for a Borel set $B\subset\Om$, 
$\capa(B)=0$ implies $\cV(B)=0$ (see \cite{Mazya64} and also \cite{Mazya8}, Sect. 12.4). 
In the present paper we will always assume  that this condition is satisfied. The operator,
associated with the closure of the form \eqref{E:form} will be denoted $H_{\cV}$.

In particular, we can consider  an  absolutely continuous measure $\cV$ which has 
a density $V\ge 0$, $V\in L^1_{loc}(\R^n)$, with
respect  to the Lebesgue measure $dx$. Such a measure will be 
absolutely continuous with respect to the capacity as well. 
 
Instead of the cubes $Q_d$ which we dealt with in Sect.\ref{S:Intro}, a more
general family of test bodies will be used. Let us start with a standard open set $\cG\subset \R^n$. 
We assume that $\cG$ satisfies the following conditions:

\ms
(a) $\cG$ is bounded and star-shaped  with respect
to an open ball $B_\rho(0)$ of radius $\rho>0$, with the center at $0\in\R^n$;

(b) $\diam(\cG)=1$.

\ms
The first condition means that $\cG$ is star-shaped with respect to every point of $B_\rho(0)$.
It implies that $\cG$ can be presented in the form
\begin{equation}\label{E:star-shaped}
\cG=\{x|\; x=r\om,\; |\om|=1,\; 0\le r< r(\om) \},
\end{equation}
where $\om\mapsto r(\om)\in (0,+\infty)$ is a Lipschitz function on the standard unit sphere 
$S^{n-1}\subset\R^n$ (see \cite{Mazya8}, Lemma 1.1.8).

The condition (b) is imposed for convenience of formulations.

For any positive $d>0$ denote by $\cG_d(0)$ the body 
$\{x|\;d^{-1}x\in\cG\}$ 
which is  homothetic to $\cG$ with coefficient $d$  and with the center of homothety at $0$. 
We will denote by $\cG_d$ a body which is obtained from $\cG_d(0)$ by a parallel
translation: $\cG_d(y)=y+\cG_d(0)$ where 
$y$ is an arbitrary  vector in $\R^n$.  

The notation $\cG_d\to\infty$ means that
the distance from $\cG_d$ to $0$ goes to infinity.

\begin{definition}\label{D:F-neg}
{\rm
Let $\ga\in(0,1)$. The {\it negligibility class} $\cN_{\ga}(\cG_d;\Om)$ consists
of all compact sets $F\subset\bcG_d$ satisfying the following conditions:
\begin{equation}\label{E:F-inclusion}
\bar\cG_d\setminus\Om\subset F\subset \bar \cG_d\;,
\end{equation}
and
\begin{equation}\label{E:G-negligibility}
\capa(F)\le \ga\capa(\bcG_d).
\end{equation}
}
\end{definition}

Now we formulate our main result about the discreteness of spectrum.

\begin{theorem}\label{T:discr}

{\rm (i)} {\rm (Necessity)} Let the spectrum of $H_\cV$ be discrete. Then for
every  function $\ga:(0,+\infty)\to (0,1)$ and every $d>0$
\begin{equation}\label{E:inf-cond}
\inf_{F\in\cN_{\ga(d)}(\cG_d,\Om)}\; \cV({\bcG_{d}\setminus F)}\to +\infty
\quad {\rm as} \quad \cG_d\to\infty.
\end{equation}

{\rm (ii)} {\rm(Sufficiency)}
Let a function $d\mapsto\ga(d)\in(0,1)$ be defined for $d>0$ 
in a neighborhood of $0$, and satisfy \eqref{E:gamma-cond1}.
Assume that 
there exists $d_0>0$
such that \eqref{E:inf-cond} holds for every $d\in (0,d_0)$.
Then the spectrum of $H_\cV$ in $L^2(\Om)$ is discrete.

\end{theorem}

Let us make some comments about this theorem.

\begin{remark}\label{R:sequence}
{\rm It suffices for the discreteness of spectrum of $H_{\cV}$
that the condition \eqref{E:inf-cond} holds 
only for a sequence of $d$'s, i.e. $d\in\{d_1,d_2,\dots\}$, $d_k\to 0$ and 
$d_k^{-2}\ga(d_k)\to +\infty$ as $k\to +\infty$.
}
\end{remark}

\begin{remark}\label{R:other-suff}
{\rm
As we will see in the proof, in the sufficiency part the condition \eqref{E:inf-cond} can be
replaced  by a weaker requirement: there exist $c>0$ and $d_0>0$  
such that for every $d\in (0,d_0)$ 
there exists $R>0$ such that 
\begin{equation}\label{E:inf-ineq}
d^{-n}\inf_{F\in\cN_{\ga(d)}(\cG_d,\Om)}\; \cV({\bcG_{d}\setminus F)}
 \ge cd^{-2}\ga(d),
\end{equation}
whenever $\bcG_d\cap(\Om\setminus B_R(0))\ne\emptyset$ 
(i.e. for distant bodies $\cG_d$ having non-empty intersection with $\Om$). 
Moreover, it suffices that the condition
\eqref{E:inf-ineq} is satisfied for a sequence $d=d_k$ satisfying the condition
formulated in Remark \ref{R:sequence}.

Note that unlike \eqref{E:inf-cond}, the condition \eqref{E:inf-ineq} does not require
that the left hand side goes to $+\infty$ as $\cG_d\to\infty$. What is actually required 
is that the left-hand side has a certain lower bound, depending on $d$ for arbitrarily small
$d>0$ and distant test bodies $\cG_d$. Nevertheless,  
the conditions \eqref{E:inf-cond} and \eqref{E:inf-ineq} 
are equivalent because each of them is equivalent to the discreteness of spectrum.
}
\end{remark}

\begin{remark}\label{R:Molchanov}
{\rm
If we take $\ga=const\in(0,1)$, then Theorem \ref{T:discr} gives Molchanov's result, but 
with the constant $\ga=c_n$
replaced by an arbitrary constant $\ga\in (0,1)$. So Theorem \ref{T:discr} contains an answer
to the above-mentioned Gelfand's question.
}
\end{remark}

\begin{remark}\label{R:ga-equiv}
{\rm
For any two functions $\ga_1,\ga_2:(0,+\infty)\to (0,1)$ satisfying the requirement 
\eqref{E:gamma-cond1}, the conditions \eqref{E:inf-cond} are equivalent, and so are 
the conditions \eqref{E:inf-ineq}, because any of these conditions is equivalent 
to the discreteness of spectrum. In a different context an equivalence of this kind
was first established in \cite{KMS}.

It follows that the conditions  \eqref{E:inf-cond} for different constants $\ga\in(0,1)$ 
are equivalent. In the particular case, when the measure $\cV$ is absolutely continuous 
with respect to the Lebesgue measure, we see that the conditions \eqref{E:main-cond}
with different constants $\ga\in(0,1)$ are equivalent. 
}
\end{remark}

\begin{remark}\label{R:domains}
{\rm The results above are new even for the 
operator $H_0=-\De$ in $L^2(\Om)$
(but for an arbitrary open set $\Om\subset \R^n$ 
with the Dirichlet boundary conditions on $\pa\Om$).
In this case the discreteness of spectrum  
is completely determined  by the geometry of $\Om$. 
Namely, for the discreteness of spectrum of $H_0$
in $L^2(\Om)$ it is 
necessary and 
sufficient that there exists $d_0>0$
such that for every $d\in (0,d_0)$
\begin{equation}\label{E:omega-cond}
\liminf_{\cG_d\to\infty}\capa(\bcG_d\setminus \Om)\ge \ga(d)\capa(\bcG_d),
\end{equation} 
where $d\mapsto \ga(d)\in (0,1)$ is a function, which is  defined in a neighborhood 
of $0$ and satisfies \eqref{E:gamma-cond1}. The conditions \eqref{E:omega-cond} 
with different functions $\ga$, satisfying the conditions above, are equivalent.
This is a non-trivial property of capacity. 
It is necessary for the discreteness of spectrum  that \eqref{E:omega-cond} 
holds for every function $\ga:(0,+\infty)\to (0,1)$ and every $d>0$, but this condition
may not be sufficient if $\ga$ does not satisfy \eqref{E:gamma-cond1}
(see Theorem \ref{T:precise} below).
}
\end{remark}

The following result demonstrates that the condition \eqref{E:gamma-cond1} is precise.

\begin{theorem}\label{T:precise}
Assume that $\ga(d)=O(d^2)$ as $d\to 0$. Then there exist an open set 
$\Om\subset\R^n$ and $d_0>0$ such that for every $d\in(0,d_0)$
the condition  \eqref{E:omega-cond} is satisfied
but the spectrum of $-\De$ in $L^2(\Om)$ with the Dirichlet boundary conditions 
is not discrete.
\end{theorem} 

Now we will state our positivity result. We will say that the operator $H_\cV$ 
is {\it strictly positive} if its spectrum does not contain $0$. 
Equivalently, we can say that the spectrum is separated from $0$. 
Since $H_{\cV}$ is defined 
by the quadratic form \eqref{E:form}, the strict positivity is equivalent
to the existence of $\la>0$ such that 
\begin{equation}\label{E:form-pos}
h_{\cV}(u,u)\ge \la \|u\|^2_{L^2(\Om)},\quad u\in C_0^\infty(\Om).
\end{equation}

\begin{theorem}\label{T:positivity}
{\rm (i) (Necessity)} Let us assume that $H_{\cV}$ is strictly positive, so that 
\eqref{E:form-pos} is satisfied with a constant $\la>0$. Let us take an arbitrary
$\ga\in (0,1)$.
Then there exist $d_0>0$ and $\vka>0$ such that 
\begin{equation}\label{E:pos-cond2}
d^{-n}\inf_{F\in \cN_{\ga}(\cG_d,\Om)} \cV(\bcG_d\setminus F)\ge \varkappa
\end{equation}
for every $d>d_0$ and every $\cG_d$.

{\rm (ii) (Sufficiency)} Assume that there exist $d>0$, $\varkappa>0$ and $\ga\in (0,1)$,
such that \eqref{E:pos-cond2} is satisfied for every $\cG_d$.
Then the operator $H_{\cV}$ is strictly positive.

Instead of all bodies $\cG_d$ it is sufficient to take only the ones from
a finite multiplicity covering (or tiling) of $\R^n$.
\end{theorem}

\begin{remark}
{\rm
Considering the Dirichlet Laplacian $H_0=-\De$ in $L^2(\Om)$
we see from Theorem \ref{T:positivity} that for any
choice of a constant $\ga\in (0,1)$ and a standard body $\cG$, 
the strict positivity of $H_0$
is equivalent to the following condition:
\begin{equation}\label{E:cap-prop}
\textrm{$\exists\,d>0$, such that 
$\capa(\bcG_d\cap(\R^n\setminus\Om))\ge\ga\capa(\bcG_d)$
for all $\cG_d$.}
\end{equation}
In particular, it follows that for two different $\ga$'s these conditions are
equivalent. Noting that $\R^n\setminus\Om$ can be an arbitrary closed subset in $\R^n$,
we get a property of the Wiener capacity, which is 
obtained as a byproduct of our spectral theory arguments.
}
\end{remark}

\section{Discreteness of spectrum: necessity}\label{S:nec}

In this section we will prove the necessity part (i) of Theorem \ref{T:discr}.
We will start by recalling some  definitions and introducing necessary notations.

For every subset $\cD\subset\R^n$ denote by $\Lip(\cD)$ the space of 
(real-valued) functions satisfying
the uniform Lipschitz condition in $\cD$, and by $\Lip_c(\cD)$ the subspace in $\Lip(\cD)$
of all functions with compact support in $\cD$ (this will be only used when $\cD$ is open).
By $\Lip_{loc}(\cD)$ we will denote the set of functions on (an open set) $\cD$ which are
Lipschitz on any compact subset $K\subset\cD$.
Note that $\Lip(\cD)=\Lip(\bar\cD)$ for any bounded $\cD$.

If $F$ is a compact subset in  an open set $\cD\subset{\R}^n$, then the 
Wiener capacity of $F$ with respect to $\cD$
is defined as
\begin{equation}\label{E:cap-def}
{\capa}_\cD (F)=\inf\left\{\left.
\int_{\R^n}|\nabla u(x)|^2 dx\,\right|\; u\in\Lip_c(\cD), u|_F=1\right\}.
\end{equation}

By $B_d(y)$ we will denote an open  ball of radius $d$ centered at $y$ in $\R^n$.
We will write $B_d$ for a ball $B_d(y)$ with unspecified center $y$.

We will  use the notation $\capa(F)$ for $\capa_{\R^n}(F)$
if $F\subset \R^n$, $n\ge 3$, and for 
$\capa_{B_{2d}}(F)$ if $F\subset \bB_d\subset \R^2$,  
where the discs $B_{d}$ and $B_{2d}$ have the same center.
The choice of these discs will be usually clear from the context,
otherwise we will specify them explicitly.

Note that the infimum does not change if we  restrict ourselves to the Lipschitz 
functions $u$ such that $0\le u\le 1$  everywhere (see e.g. \cite{Mazya8}, Sect. 2.2.1).

\ms
We will also need another (equivalent) definition of the Wiener capacity $\capa(F)$ 
for a compact set $F\subset \bB_d$. For  $n\ge 3$ it is as follows:
\begin{equation}\label{E:cap-def2}
\capa(F)=\sup\{\mu(F)\left| \int_F\cE(x-y) d\mu(y)\le 1\quad \text{on}\ \R^n\setminus F\right.\},
\end{equation}
where the supremum is taken over all positive finite Radon measures $\mu$ on $F$ and 
$\cE=\cE_n$ is the standard fundamental solution of
$-\De$ in
$\R^n$ i.e.
\begin{equation}\label{E:fund-sol}
\cE(x)=\frac{1}{(n-2)\om_n}|x|^{2-n}\;,
\end{equation}
with $\om_n$ being the area of the unit sphere $S^{n-1}\subset\R^n$. If $n=2$, then
\begin{equation}\label{E:cap-def2-2}
\capa(F)=\sup\{\mu(F)\left| \int_F G(x,y) d\mu(y)\le 1\quad \text{on}\ B_{2d}\setminus F\right.\},
\end{equation}
where $G$ is the Green function of the Dirichlet problem for $-\De$ in $B_{2d}$, i.e. 
\begin{equation*}\label{E:Green-function-2}
-\De G(\cdot -y)=\de(\cdot -y), \quad y\in B_{2d}, 
\end{equation*}
$G(\cdot,y)|_{\pa B_{2d}}=0$ for all $y\in B_{2d}$.
The maximizing measure in \eqref{E:cap-def2} or in \eqref{E:cap-def2-2} exists and is unique.
We will denote it $\mu_F$ and call it the \textit{equilibrium measure}.
Note that
\begin{equation*}\label{E:mu-F-property}
\capa(F)=\mu_F(F)=\mu_F(\R^n)=\langle\mu_F,1\rangle.
\end{equation*}
The corresponding potential will be denoted $P_F$, so 
\begin{equation*}\label{E:psi-F}
P_F(x)=\int_F\cE(x-y)d\mu_F(y), \quad x\in \R^n\setminus F, \qquad n\ge 3,
\end{equation*} 
\begin{equation*}\label{E:psi-F-2}
P_F(x)=\int_F G(x,y) d\mu_F(y), \quad x\in B_{2d}\setminus F, \qquad n=2.
\end{equation*}
We will call $P_F$ the \textit{equilibrium potential} or \textit{capacitary potential}. 
We will  extend it to $F$ by setting $P_F(x)=1$ for all $x\in F$.

It follows from the maximum principle that $0\le P_F\le 1$
everywhere in $\R^n$ if $n\ge 3$ (and in $B_{2d}$ if $n=2$).

In case when  $F$ is a closure
of an open subset with a smooth boundary, $u=P_F$ is the unique minimizer
for the Dirichlet integral in \eqref{E:cap-def} where we should take
$\cD=\R^n$ if $n\ge 3$ and $\cD=B_{2d}$ if $n=2$. In particular, 
\begin{equation}\label{E:minimizer}
\int |\nabla P_F|^2 dx=\capa(F),
\end{equation}
where the integration is taken over $\R^n$ (or $\R^n\setminus F$) if $n\ge 3$ and
over $B_{2d}$ (or $B_{2d}\setminus F$) if $n=2$.

\medskip
The following lemma provides an  auxiliary estimate
which is needed for the proof.
 
\begin{lemma}\label{L:grad-P-est} 
Assume that $\cG$ has a $C^\infty$ boundary, and $P$ is the equilibrium potential of $\bcG_d$.
Then
\begin{equation}\label{E:grad-P-est}
\int_{\pa\cG_d}|\nabla P|^2 ds \le nL\rho^{-1}d^{-1}\capa(\bar \cG_d),
\end{equation}
where  the gradient $\nabla P$ in the left hand side is taken 
along the exterior of $\bcG_d$, $ds$ is the $(n-1)$-dimensional volume element
on $\pa\cG_d$.
The positive constants $\rho, L$ are geometric characteristics of the standard body $\cG$
(they depend on the choice of 
$\cG$ only, but not on $d$):
$\rho$ was introduced at the beginning of Section \ref{S:main}, and
\begin{equation}\label{E:L}
L=\left[\inf_{x\in\pa\cG} \nu_r(x)\right]^{-1},
\end{equation}
where $\nu_r(x)=\frac{x}{|x|}\cdot \nu(x)$, $\nu(x)$ is the unit normal vector to 
$\pa\cG$ at $x$ which is directed to the exterior of $\bcG$.
\end{lemma}

{\bf Proof.} It suffices to consider $\cG_d=\cG_d(0)$. For simplicity we will write $\cG$ 
instead of $\cG_d(0)$ in this proof, until the size becomes relevant. 

We will first consider the case $n\ge 3$. Note that $\De P=0$ on 
$\cbcG=\R^n\setminus\bcG$. Also $P=1$ on $\bcG$, so in fact $|\nabla P|=|\pa P/\pa\nu|$.
Using the Green formula, we obtain
\begin{align*}
0 & = \int_{\cbcG}\De P\cdot\frac{\pa P}{\pa r} dx
= \int_{\cbcG} \De P \left(\frac{x}{|x|}\cdot\nabla P\right) dx\\
& = -\int_{\cbcG} \nabla P\cdot \nabla\left(\frac{x}{|x|}\cdot\nabla P\right) dx
- \int_{\pa\cG}\frac{\pa P}{\pa\nu}\left(\frac{x}{|x|}\cdot\nabla P\right)ds\\
& = -\sum_{i,j}\int_{\cbcG}\frac{\pa P}{\pa x_j}\cdot\frac{\pa}{\pa x_j}
\left(\frac{x_i}{|x|}\cdot\frac{\pa P}{\pa x_i}\right) dx 
- \int_{\pa\cG}\frac{\pa P}{\pa\nu}\cdot\frac{\pa P}{\pa r}ds\\
& = -\sum_{i,j}\int_{\cbcG}\frac{\pa P}{\pa x_j}\cdot
\frac{\de_{ij}}{|x|}\cdot\frac{\pa P}{\pa x_i} dx + 
\sum_{i,j}\int_{\cbcG}
\frac{x_i x_j}{|x|^3}\cdot\frac{\pa P}{\pa x_i} \cdot \frac{\pa P}{\pa x_j} dx\\
& -\sum_{i,j}\int_{\cbcG}
\frac{x_i}{|x|}\cdot\frac{\pa P}{\pa x_j} \cdot \frac{\pa^2 P}{\pa x_i\pa x_j} dx
- \int_{\pa \cG}\frac{\pa P}{\pa\nu}\cdot\frac{\pa P}{\pa r}ds\\
&= -\int_{\cbcG}\frac{1}{|x|}|\nabla P|^2 dx + 
\int_{\cbcG}\frac{1}{|x|}\left|\frac{\pa P}{\pa r}\right|^2 dx\\
&-\frac{1}{2}\sum_{i}\int_{\cbcG}
\frac{x_i}{|x|}\cdot\frac{\pa }{\pa x_i}|\nabla P|^2 dx
- \int_{\pa \cG}|\nabla P|^2\nu_r ds.
\end{align*}
Integrating by parts in the last integral over $\cbcG$, 
we see that it equals
\begin{align*}
&\frac{1}{2}\sum_{i}\int_{\cbcG}
\frac{\pa }{\pa x_i}\left(\frac{x_i}{|x|}\right)\cdot|\nabla P|^2 dx 
+\frac{1}{2}\sum_i\int_{\pa \cG}\frac{x_i}{|x|}|\nabla P|^2\nu_i ds\\
&=\frac{n-1}{2}\int_{\cbcG}\frac{1}{|x|}|\nabla P|^2dx
+\frac{1}{2}\int_{\pa \cG}|\nabla P|^2\nu_r ds,
\end{align*}
where $\nu_i$ is the $i$th component of $\nu$. Returning to the 
calculation above, we obtain
\begin{equation}\label{E:nabla-P-identity}
0=\frac{n-3}{2}\int_{\cbcG}\frac{1}{|x|}|\nabla P|^2 dx + 
\int_{\cbcG}\frac{1}{|x|}\left|\frac{\pa P}{\pa r}\right|^2 dx
-\frac{1}{2}\int_{\pa \cG}|\nabla P|^2\nu_r ds.
\end{equation}
It follows that
\begin{equation*}\label{E:nabla-P-ineq}
\int_{\pa \cG}|\nabla P|^2 \nu_r ds \le (n-1)\int_{\cbcG}\frac{1}{|x|}|\nabla P|^2 dx.
\end{equation*}
Recalling that $\cG=\cG_d(0)$, we observe that $|x|^{-1}\le (\rho d)^{-1}$. Now using  
\eqref{E:minimizer}, we obtain the desired estimate \eqref{E:grad-P-est}  for $n\ge 3$
(with $n-1$ instead of $n$).

Let us consider the case $n=2$. Then, by definition, the equilibrium potential  $P$ for $\cG=\cG_d(0)$
is defined in the ball $B_{2d}(0)$. It satisfies $\De P=0$ in $B_{2d}(0)\setminus \bcG$ and the boundary
conditions $P|_{\pa\cG}=1$,
$P|_{\pa B_{2d}(0)}=0$. Let us first modify the calculations above by taking the integrals over
$B_{\de}(0)\setminus \bcG$ (instead of $\cbcG$), where
$d<\de<2d$.
We will get additional boundary terms with the integration
over $\pa B_{\de}(0)$. Instead of \eqref{E:nabla-P-identity} we will obtain
\begin{align*}\label{E:nabla-P-identity2}
&0=-\frac{1}{2}\int_{B_{\de}(0)\setminus \bcG}\frac{1}{|x|}|\nabla P|^2 dx + 
\int_{B_{\de}(0)\setminus \bcG}\frac{1}{|x|}\left|\frac{\pa P}{\pa r}\right|^2 dx\\
&-\frac{1}{2}\int_{\pa \cG}|\nabla P|^2\nu_r ds 
+ \frac{1}{2}\int_{\pa B_{\de}(0)}\left[2\left|\frac{\pa P}{\pa r}\right|^2-|\nabla P|^2\right] ds.
\end{align*}
Therefore
\begin{align*}
\int_{\pa \cG}|\nabla P|^2\nu_r ds 
&\le \int_{B_{\de}(0)\setminus \bcG}\frac{1}{|x|}|\nabla P|^2 dx +
\int_{\pa B_{\de}(0)}\left[2\left|\frac{\pa P}{\pa r}\right|^2-|\nabla P|^2\right] ds\\ 
&\le \frac{1}{\rho d}\int_{B_{2d}(0)\setminus \bcG}|\nabla P|^2 dx +
\int_{\pa B_{\de}(0)}|\nabla P|^2 ds.
\end{align*}
Now let us integrate both sides with respect to $\de$ over the interval 
$[d,2d]$ and divide the result by $d$ (i.e. take average over all $\de$). 
Then the left hand side and the first term in the right hand side do not change, 
while the last term becomes $d^{-1}$ times the volume integral with respect  
to the Lebesgue measure over $B_{2d}(0)\setminus B_d(0)$. Due to 
\eqref{E:minimizer} the right hand side can be estimated by 
$(1+\rho)(\rho d)^{-1}\capa(\bcG_d)$. Since $0<\rho\le 1$, we get the estimate
\eqref{E:grad-P-est} for $n=2$. $\square$

\bs
{\bf Proof of Theorem \ref{T:discr}, part (i).} (a) We will use the same notations as above.
Let us fix $d>0$, take $\cG_d=\cG_d(z)$, and assume that $\cG$ has a $C^\infty$ boundary. 
Let us take a compact set $F\subset \R^n$ with the following properties:

\ms
(i) $F$ is the closure  of an open set with a $C^\infty$ boundary;

(ii) $\bcG_d\setminus \Om\Subset F\subset B_{3d/2}(z)$;

(iii) $\capa(F)\le \ga\capa(\bcG_d)$ with $0<\ga<1$. 

\msni
Let us recall that the notation $\bcG_d\setminus \Om\Subset F$ means that $\bcG_d\setminus \Om$
is contained in the interior of $F$. This implies that $\cV(\bcG_d\setminus F)<+\infty$.
The inclusion $F\subset B_{3d/2}(z)$ and the inequality (iii) hold, in particular, 
for compact sets $F$ which are
small neighborhoods  (with smooth boundaries) of negligible compact subsets of $\bcG_d$, 
and it is exactly such $F$'s which we have in mind.

We will refer to the sets $F$ satisfying (i)-(iii) above as \textit{regular} ones.

Let $P$ and $P_F$ denote the
equilibrium potentials  of $\bcG_d$ and $F$ respectively. The equilibrium measure 
$\mu_{\bcG_d}$ has its support in 
$\pa\cG_d$ and has density $-\pa P/\pa\nu$ with respect to the $(n-1)$-dimensional Riemannian
measure $ds$ on $\pa\cG_d$. So for $n\ge 3$ we have
\begin{equation*}
P(y)= -\int_{\pa\cG_d}\cE(x-y) \frac{\pa P}{\pa\nu}(x) ds_x, \quad y\in \R^n;
\end{equation*}
\begin{equation*}
-\int_{\pa\cG_d}\frac{\pa P}{\pa\nu}(x) ds_x = \capa(\bcG_d);
\end{equation*}
\begin{equation*}
P(y) = 1 \ {\rm for\ all} \ y\in \cG_d,\quad 0\le P(y)\le 1 \ {\rm for \ all}\ y\in\R^n.
\end{equation*}
(If $n=2$, then the same holds only with $y\in B_{2d}$ and with the fundamental solution 
$\cE$ replaced by the Green function
$G$.) It follows that
\begin{equation*}
-\int_{\pa\cG_d} P_F \frac{\pa P}{\pa\nu}ds=
-\int_F \int_{\pa\cG_d}\cE(x-y)\frac{\pa P}{\pa\nu}(x)ds_x d\mu_F(y)\le \mu_F(F)=\capa(F).
\end{equation*}
Therefore, 
\begin{equation*}
\capa(\bcG_d)-\capa(F)\le -\int_{\pa\cG_d}(1-P_F)\frac{\pa P}{\pa\nu} ds,
\end{equation*}
and, using Lemma \ref{L:grad-P-est}, we obtain
\begin{align}\label{E:cap-cap}
&(\capa(\bcG_d)-\capa(F))^2
\le \left(\int_{\pa\cG_d}(1-P_F)\frac{\pa P}{\pa\nu} ds\right)^2\\
& \le\|1-P_F\|^2_{L^2(\pa\cG_d)}\|\nabla P\|^2_{L^2(\pa\cG_d)}\le 
nL(\rho d)^{-1}\capa(\cG_d)\|1-P_F\|^2_{L^2(\pa\cG_d)}, \notag
\end{align}
where $L$ is defined by \eqref{E:L}.

\ms
(b) Our next goal will be to estimate the norm $\|1-P_F\|_{L^2(\pa\cG_d)}$ in \eqref{E:cap-cap} by 
the norm of the same function in $L^2(\cG_d)$. 
We will use the polar coordinates
$(r,\om)$ as in
\eqref{E:star-shaped},  so in particular $\pa\cG_d$ is presented as the set $\{r(\om)\om|\;\om\in S^{n-1}\}$,
where $r:S^{n-1}\to (0,+\infty)$ is a Lipschitz function 
($C^\infty$ as long as we assume  the boundary $\pa\cG$ to be $C^\infty$). 
Assuming that $v\in \Lip(\bcG_d)$,
we can write
\begin{align}\label{E:int-v2}
\int_{\pa\cG_d}|v|^2 ds &= \int_{S^{n-1}} |v|^2 \frac{r(\om)^{n-1}}{\nu_r}d\om\\
&\le L\int_{S^{n-1}}|v(r(\om),\om)|^2 r(\om)^{n-1}d\om, \notag
\end{align}
where $d\om$ is the standard $(n-1)$-dimensional volume element on $S^{n-1}$.

Using the inequality
\begin{equation*}
|f(\eps)|^2\le 2\eps\int_0^\eps |f'(t)|^2 dt + \frac{2}{\eps}\int_0^\eps |f(t)|^2 dt, 
\quad f\in \Lip([0,\eps]),\quad \eps>0,
\end{equation*}
we obtain
\begin{align*}
&|v(r(\om),\om)|^2\\
&\le 2\eps r(\om)\int_{(1-\eps)r(\om)}^{r(\om)}|v'_{\rho}(\rho,\om)|^2  d\rho
+\frac{2}{\eps r(\om)}\int_{(1-\eps)r(\om)}^{r(\om)} |v(\rho,\om)|^2 d\rho \\
&\le \frac{2\eps r(\om)}{[(1-\eps)r(\om)]^{n-1}}
\int_{(1-\eps)r(\om)}^{r(\om)}|v'_{\rho}(\rho,\om)|^2 \rho^{n-1} d\rho\\
& +\frac{2}{\eps r(\om)[(1-\eps)r(\om)]^{n-1}}\int_{(1-\eps)r(\om)}^{r(\om)} |v(\rho,\om)|^2 \rho^{n-1} d\rho. 
\end{align*}
It follows that the integral in the right hand side of \eqref{E:int-v2}
is estimated by
\begin{align*}
&\int_{S^{n-1}}\frac{2\eps r(\om)d\om}{(1-\eps)^{n-1}}
\int_{(1-\eps)r(\om)}^{r(\om)}|v'_{\rho}(\rho,\om)|^2 \rho^{n-1} d\rho\\
&+\int_{S^{n-1}}\frac{2 d\om}{\eps (1-\eps)^{n-1}r(\om)} |v(\rho,\om)|^2 \rho^{n-1} d\rho.
\end{align*}
Taking $\eps \le 1/2$, we can majorize this by
\begin{align*}
2^n\eps d\int_{\bcG_d}|\nabla v|^2 dx + \frac{2^n}{\eps\rho d}\int_{\bcG_d}|v|^2 dx,
\end{align*}
where $\rho\in (0,1]$ is the constant from the description of $\cG$ in Sect. \ref{S:main}.
Recalling \eqref{E:int-v2}, we see that the resulting estimate
has the form
\begin{equation*}\label{E:int-v2-final}
\int_{\pa\cG_d}|v|^2 ds \le 
2^n L \eps d\int_{\bcG_d}|\nabla v|^2 dx + \frac{2^n L}{\eps\rho d}\int_{\bcG_d}|v|^2 dx.
\end{equation*} 
Now, taking $v=1-P_F$, 
we obtain   
\begin{equation*}\label{E:int-PF2}
\int_{\pa\cG_d}(1-P_F)^2 ds \le 
2^n L \eps d\capa(F) + \frac{2^n L}{\eps\rho d}\int_{\bcG_d}(1-P_F)^2 dx.
\end{equation*}
Using this estimate in \eqref{E:cap-cap}, we obtain
\begin{align}\label{E:cap-cap2}
&(\capa(\bcG_d)-\capa(F))^2\\
&\le \rho^{-1}n2^n L^2\capa(\bcG_d)\left(\eps \capa(F)+
\frac{1}{\eps\rho d^2}\int_{\cG_d}(1-P_F)^2 dx\right). \notag
\end{align}

\ms
(c) Now let us consider $\cG$ which is star-shaped
with respect to a ball, but not necessarily has $C^\infty$ boundary.
In this case we can approximate the function $r(\om)$ (see Section \ref{S:main}) from above by 
a decreasing sequence of $C^\infty$ functions $r_k(\om)$ (e.g. we can apply a standard mollifying
procedure to $r(\om)+1/k$), so that for the the corresponding bodies $\cG^{(k)}$
the constants $L_k$ are uniformly bounded. 
It is clear that in this case we will also have $\rho_k\ge \rho$, 
and $\capa(\bcG_d^{(k)})\to\capa(\bcG_d)$ due to the well known continuity property of the capacity
(see e.g. Section 2.2.1 in \cite{Mazya8}). So we can pass to the limit in \eqref{E:cap-cap2}
as $k\to +\infty$ and conclude that it holds for arbitrary $\cG$ (which is star-shaped with respect 
to a ball).  But for the moment we still retain the regularity condition on $F$.

\ms
(d) Let us define
\begin{equation}\label{E:cL}
\cL=\left\{u\left|u\in C_0^\infty(\Om),\; h_\cV(u,u)+\|u\|^2_{L^2(\Om)} 
\le 1\right.\right\},
\end{equation}
where $h_\cV$ is defined by \eqref{E:form}.
By the standard functional analysis argument (see e.g. Lemma 2.3 in \cite{Kondrat'ev-Shubin}) 
the spectrum  of $H_{\cV}$ is discrete if and only if $\cL$ is precompact in $L^2(\Om)$,
which in turn holds if and only if $\cL$ has ``small tails", 
i.e. for every $\eta>0$ there exists $R>0$ such that 
\begin{equation}\label{E:smalltails}
\int_{\Om\setminus B_R(0)}|u|^2 dx\le \eta \quad \text{for every} \quad u\in\cL,
\end{equation}
Equivalently, we can write that
\begin{equation}\label{E:smalltails2}
\int_{\Om\setminus B_R(0)}|u|^2 dx\le \eta \left[\int_{\Om}|\nabla u|^2 dx +
\int_{\Om} |u|^2\cV(dx)\right],
\end{equation} 
for every  $u\in C_0^\infty(\Om)$. 

Therefore, it follows from the discreteness of spectrum of $H_\cV$ that for every $\eta>0$ there exists
$R>0$ such that for every $\cG_d$ with $\bcG_d\cap (\R^n\setminus B_R(0))\ne\emptyset$
and every $u\in C_0^\infty(\cG_d\cap\Om)$  
\begin{equation}\label{E:mu-estimate}
\int_{\cG_d}|u|^2 dx\le \eta\left(\int_{\cG_d}|\nabla u|^2 dx + 
\int_{\bcG_d}|u|^2\cV(dx)\right). 
\end{equation}
In other words, $\eta=\eta(\cG_d)\to 0$ as $\cG_d\to \infty$  for the best constant in 
\eqref{E:mu-estimate}. (Note that $\eta(\cG_d)^{-1}$ is the bottom of the Dirichlet spectrum of
$H_\cV$ in $\cG_d\cap\Om$.)

Since  $1-P_F=0$ on $F$ (hence in a neighborhood of $\bcG_d\setminus \Om$),
we can take $u=\chi_\sigma(1-P_F)$, where $\sigma\in (0,1)$ to be chosen later,
$\chi_\sigma\in C_0^\infty(\cG_d)$ is a cut-off function satisfying
$0\le \chi_\sigma\le 1$, $\chi_\sigma=1$ on $\cG_{(1-\sigma)d}$, and $|\nabla\chi_\sigma|\le Cd^{-1}$
with $C=C(\cG)$. Then, using integration by parts and the equation $\De P_F=0$  on $\cG\setminus F$,
we obtain
\begin{align*}
\int_{\cG_d}|\nabla u|^2 dx  
&=\int_{\cG_d}\left(|\nabla\chi_\sigma|^2(1-P_F)^2
-\nabla(\chi_\sigma^2)\cdot (1-P_F)\nabla P_F
+ \chi_\sigma^2|\nabla P_F|^2  \right) dx\\
&=\int_{\cG_d}|\nabla\chi_\sigma|^2(1-P_F)^2 dx 
\le C^2(\sigma d)^{-2}\int_{\cG_d} (1-P_F)^2 dx.
\end{align*}
Therefore, from \eqref{E:mu-estimate}  
\begin{equation*}
\int_{\cG_d}|u|^2 dx
\le \eta\left[C^2(\sigma d)^{-2}\int_{\cG_d} (1-P_F)^2 dx + 
\cV(\bcG_d\setminus F)\right],
\end{equation*} 
hence
\begin{equation*}
\int_{\cG_{(1-\sigma)d}}(1-P_F)^2 dx
\le \eta\left[C^2(\sigma d)^{-2}\int_{\cG_d} (1-P_F)^2 dx + 
\cV(\bcG_d\setminus F)\right].
\end{equation*} 
Now, applying the obvious estimate
\begin{align*}
\int_{\cG_{d}}(1-P_F)^2 dx&\le \int_{\cG_{(1-\sigma)d}}(1-P_F)^2 dx 
+\mes(\cG_d\setminus\cG_{(1-\sigma)d})\\
&\le\int_{\cG_{(1-\sigma)d}}(1-P_F)^2 dx + C_1\sigma d^n,
\end{align*} 
with $C_1=C_1(\cG)$, we see that
\begin{equation*}
\int_{\cG_{d}}(1-P_F)^2 dx\le \eta\left[C^2(\sigma d)^{-2}\int_{\bcG_d} (1-P_F)^2 dx + 
\cV(\bcG_d\setminus F)\right]+C_1\sigma d^n,
\end{equation*}
hence
\begin{equation}\label{E:PF-eta}
\int_{\cG_{d}}(1-P_F)^2 dx\le 2\eta \cV(\bcG_d\setminus F) + 2C_1\sigma d^n,
\end{equation}
provided
\begin{equation}\label{E:eta-C2}
\eta C^2 (\sigma d)^{-2}\le 1/2.
\end{equation}
Returning to \eqref{E:cap-cap2} and using \eqref{E:PF-eta} we obtain
\begin{align}\label{E:1-cap/cap}
\left(1-\frac{\capa(F)}{\capa(\bcG_d)}\right)^2
\le C_2\left[\eps+\eps^{-1}d^{-n}\int_{\cG_d}(1-P_F)^2 dx\right]\\
\le C_2[\eps + 2C_1\sigma\eps^{-1} + 2\eps^{-1} d^{-n}\eta\cV(\bcG_d\setminus F)],
\notag
\end{align}
where $C_2=C_2(\cG)$. Without loss of generality we will assume that $C_2\ge 1/2$.
Recalling that $\capa(F)\le \ga \capa(\bcG_d)$, we can replace the ratio
$\capa(F)/\capa(\bcG_d)$ in the left hand side by $\ga$. Now let us choose
\begin{equation}\label{E:epsilon}
\eps=\frac{(1-\ga)^2}{4C_2},\quad \sigma=\frac{\eps(1-\ga)^2}{8C_1}=\frac{(1-\ga)^4}{32 C_1 C_2}.
\end{equation}
Then $\eps\le 1/2$ and for every fixed $\gamma\in (0,1)$ and $d>0$ 
the condition \eqref{E:eta-C2} will be satisfied for distant bodies $\cG_d$,
because $\eta=\eta(\cG_d)\to 0$ as $\cG_d\to\infty$.
(More precisely, there exists $R=R(\gamma, d)>0$, such that \eqref{E:eta-C2} holds for every $\cG_d$
such that $\cG_d\cap(\R^n\setminus B_R(0))\ne\emptyset$.) 

If $\eps$ and $\sigma$ are chosen according to \eqref{E:epsilon}, then \eqref{E:1-cap/cap} becomes
\begin{equation}\label{E:main-nec-ineq}
d^{-n}\cV(\bcG_d\setminus F)\ge (16 C_2\eta)^{-1} (1-\ga)^4,
\end{equation}
which holds for distant bodies $\cG_d$ if  $\gamma\in (0,1)$ and $d>0$
are arbitrarily fixed.

\ms
(e) Up to this moment we worked with ``regular" sets $F$ --
see conditions (i)-(iii) in the part (a) of this proof. 
Now we can get rid of the regularity requirements (i) and (ii),
retaining (iii). So let us assume that $F$ is a compact set,
$\bcG_d\setminus\Om\subset F\subset \bcG_d$ and $\capa(F)\le\ga\capa(\bcG_d)$ 
with $\ga\in(0,1)$.  Let us construct 
a sequence of compact sets $F_k\Supset F$, $k=1,2,\dots$, such that 
every $F_k$ is regular, 
\begin{equation*}
F_1\Supset F_2\Supset\dots, \quad  \text{and} \quad\bigcap_{k=1}^\infty F_k=F.
\end{equation*}
We have then $\capa(F_k)\to\capa(F)$ as $k\to +\infty$ due to the 
well known continuity property  of the capacity (see e.g. Section 2.2.1 in \cite{Mazya8}). 
According to the previous steps of this proof, the inequality \eqref{E:main-nec-ineq} 
holds for distant $\cG_d$'s if we replace $F$ by $F_k$ and $\ga$ by $\ga_k=\capa(F_k)/\capa(\bcG_d)$. Since
the measure $\cV$ is positive, the resulting inequality will still hold if
we  replace $\cV(\bcG_d\setminus F_k)$ by $\cV(\bcG_d\setminus F)$. Taking limit as
$k\to +\infty$, we obtain that \eqref{E:main-nec-ineq} holds with $\ga'=\capa(F)/\capa(\bcG_d)$
instead of $\ga$. Since $\ga'\le\ga$, \eqref{E:main-nec-ineq} immediately follows for
arbitrary compact $F$ such that $\bcG_d\setminus\Om\subset F\subset \bcG_d$ and 
$\capa(F)\le\ga\capa(\bcG_d)$  with $\ga\in(0,1)$.

\ms
(f) Let us  fix $\cG$ and take infimum over all 
negligible $F$'s (i.e. compact sets $F$, such that
$\bcG_d\setminus \Om\subset F\subset\bcG_d$ and $\capa(F)\le\ga\capa(\bcG_d)$) 
in the right hand side of \eqref{E:main-nec-ineq}. We get then for distant $\cG_d$'s
\begin{equation}\label{E:main-nec-ineq-inf}
d^{-n}\inf_{F\in\cN_\ga(\cG_d,\Om)}\cV(\bcG_d\setminus F)\ge (16 C_2\eta)^{-1} (1-\ga)^4.
\end{equation} 
Now let us recall that the discreteness of spectrum is equivalent to
the condition $\eta=\eta(\cG_d)\to 0$ as $\cG_d\to\infty$ (with any fixed $d>0$). If this
is the case, then it is clear  from \eqref{E:main-nec-ineq-inf}, that 
for every fixed $\ga\in(0,1)$ and $d>0$, the left hand side of \eqref{E:main-nec-ineq-inf} 
tends to $+\infty$  as $\cG_d\to\infty$. This concludes the proof of part (i) 
of Theorem \ref{T:discr}. $\square$

\section{Discreteness of spectrum: sufficiency}\label{S:sufficient}

In this section we will establish the sufficiency part of Theorem \ref{T:discr}.

\ms
Let us recall the {\it Poincar\'e inequality}
(see e.g. \cite{Gilbarg-Trudinger}, Sect. 7.8, or \cite{Kondrat'ev-Shubin}, Lemma 5.1):
\begin{equation*}\label{E:Poincare}
||u -\bar u||^2_{L^2(\cG_d)}
\le A(\cG){d^2} \int_{\cG_d} |\nabla u(x)|^2 dx, \quad u\in\Lip(\cG_d),
\end{equation*}
where $\cG_d\subset\R^n$ was described in Section \ref{S:main}
\begin{equation*}\label{E:mean-def}
\bar u = \frac{1}{|\cG_d|} \int_{\cG_d} u(x)\, dx
\end{equation*}
is the mean value of $u$ on $\cG_d$, $|\cG_d|$ is the Lebesgue volume of $\cG_d$, 
$A(\cG)>0$ is independent of
$d$.  (In fact, the best $A(\cG)$ is obtained if  $A(\cG)^{-1}$ is the lowest
non-zero Neumann eigenvalue of  $-\Delta$ in $\cG$.) 

\ms
The following Lemma generalizes (to an arbitrary body $\cG$) a particular case  
of the first part of Theorem 10.1.2 in \cite{Mazya8} (see also Lemma 2.1 in \cite{KMS}).  

\begin{lemma}\label{L:Mazya1}
There exists $C(\cG)>0$ such that the following inequality holds for every 
function $u\in \Lip(\bcG_d)$  which
vanishes on a compact set $F\subset \bcG_d$ (but is not identically $0$ on $\bcG_d$):
\begin{equation}\label{E:cap-above}
\capa(F)\le \frac{C(\cG)\int_{\cG_d} |\nabla u(x)|^2dx}{|\cG_d|^{-1}\int_{\cG_d} |u(x)|^2 dx}\;.
\end{equation}
\end{lemma}

\medskip
{\bf Proof.}
Let us normalize $u$ by 
\begin{equation*}\label{E:normalization}
|\cG_d|^{-1}\int_{\cG_d} |u(x)|^2 dx=1,
\end{equation*}
i.e. $\overline {|u|^2}=1$.  
By the Cauchy-Schwarz inequality we obtain
\begin{equation}\label{E:|u|}
\overline{|u|}\le \left(\overline{|u|^2}\right)^{1/2}=1
\end{equation} 

Replacing $u$ by $|u|$ does not change the denominator
and may only decrease the numerator in \eqref{E:cap-above}. Therefore
we can restrict ourselves to Lipschitz functions $u\ge 0$.  

Let us denote $\phi=1-u$. Then $\phi=1$ on $F$, and $\bar\phi=1-\bar u\ge 0$ 
due to \eqref{E:|u|}. Let us estimate $\bar\phi$ from above. 
Obviously
\begin{equation*}\label{E:1-|u|}
\bar\phi=|\cG_d|^{-1/2}(\|u\|-\|\bar u\|)
\le |\cG_d|^{-1/2}\|u-\bar u\|,
\end{equation*}
where $\|\cdot\|$ means the norm in $L^2(\cG_d)$.
Hence  
the Poincar\'e inequality gives
\begin{equation*}\label{E:1-bar-u}
\bar\phi
\le A^{1/2} d|\cG_d|^{-1/2}\|\nabla u\|=A^{1/2} d|\cG_d|^{-1/2}\|\nabla \phi\|, 
\end{equation*}
where $A=A(\cG)$.
So
\begin{equation*}\label{E:phi-square}
\bar{\phi}^2
\le A d^{2} |\cG_d|^{-1} \int_{\cG_d} |\nabla\phi|^2dx.
\end{equation*}
and
\begin{equation*}
\|\bar\phi\|^2\le A d^2 \int_{\cG_d} |\nabla\phi|^2dx.
\end{equation*}

Using the Poincar\'e inequality again, we obtain
\begin{equation*}\label{E:phi-square-est}
\|\phi\|^2=\|(\phi-\bar\phi)+\bar\phi\|^2\le
2\|\phi-\bar\phi\|^2+2\|\bar\phi\|^2\le
4Ad^2 \int_{\cG_d} |\nabla\phi |^2dx,
\end{equation*}
or
\begin{equation}\label{E:Mazya(1)}
\int_{\cG_d} \phi^2 dx\le 4 A d^2 \int_{\cG_d} |\nabla \phi |^2 dx. 
\end{equation}
Let us extend $\phi$ outside $\cG_d=\cG_d(y)$ by inversion in each ray
emanating from $y$. In notations introduced in \eqref{E:star-shaped}
we can write that $\phi(y+r\om)=\phi(y+r^{-1}(r(\om))^2\om)$
for every $r>r(\om)$ and every $\om\in S^{n-1}$.

It is easy to see that the extension $\tilde \phi$ satisfies
\begin{equation*}\label{E:Q3d}
\int_{B_{3d}}|\tilde\phi|^2dx\le C_1(\cG) \int_{\cG_d}|\phi|^2dx,\quad
\int_{B_{3d}} |\nabla \tilde\phi|^2 dx\le C_1(\cG)\int_{\cG_d} |\nabla \phi|^2dx.  
\end{equation*}
Let $\eta$ be a  piecewise smooth function, such that
$\eta=1$ on $B_d$, $\eta=0$ outside $B_{2d}$,
$0\le \eta\le 1$ and $|\nabla\eta|\le d^{-1}$, i.e. 
$\eta(x)=2-d^{-1}|x|$ if $d\le |x|\le 2d$. 
Then
\begin{equation*}\label{E:capF-est}
\capa(F)\le 
\int_{B_{2d}}|\nabla (\tilde\phi\eta)|^2dx
\le 2C_1(\cG)\left(\int_{\cG_d}|\nabla \phi|^2dx+
d^{-2}\int_{\cG_d} \phi^2 dx\right).
\end{equation*}
Taking into account that $|\nabla\phi|=|\nabla u|$
and using \eqref{E:Mazya(1)}, we obtain
\begin{equation*}\label{E:capF-est2}
\capa(F)\le 2C_1(\cG)(1+4A) \int_{\cG_d} |\nabla u|^2 dx,
\end{equation*}
which is equivalent to \eqref{E:cap-above} with
$C(G)=2C_1(\cG)(1+4A(\cG))$. $\square$

\ms
The next lemma is an adaptation of a very general  Lemma 12.1.1 from \cite{Mazya8}
(see also Lemma 2.2 in \cite{KMS}) to general test bodies $\cG_d$ (instead of cubes $Q_d$). 

\begin{lemma}\label{L:Mazya2}
Let $\cV$ be a positive Radon measure in $\Om$.
There exists $C_2(\cG)>0$ such that for every $\gamma\in(0,1)$ and $u\in \Lip(\bcG_d)$ with
$u=0$ in a neighborhood of $\bcG_d\setminus \Om$,
\begin{equation}\label{E:Mazya(2)}
\int_{\cG_d}|u|^2 dx\le {C_2(\cG) d^2\over \gamma}  \int_{\cG_d}|\nabla u|^2 dx
+{C_2(\cG)d^n \over \cV_\ga(\cG_d,\Om)}\
\int_{\bcG_d}|u|^2 \cV(dx),
\end{equation}
where
\begin{equation}\label{E:functional}
\cV_{\ga}(\cG_d,\Om)=\inf_{F\in\cN_\ga(\cG_d,\Om)}\cV(\cG_d\setminus F).
\end{equation}
(Here the negligibility class $\cN_\ga(\cG_d,\Om)$ was introduced in Definition \ref{D:F-neg}.)
\end{lemma}

{\bf Proof.} Let ${\cal M}_\tau=\{x\in \bcG_d:|u(x)|> \tau\},$
where $\tau\ge 0$. Note that ${\cal M}_\tau$ is a relatively open subset of $\bcG$,
and ${\cal M}_\tau\subset \Om$, hence $\bcG_d\setminus {\cal M}_\tau\supset \bcG_d\setminus \Om$. 
Since
\begin{equation*}
|u|^2\le 2\tau^2+2(|u|-\tau)^2 \quad \hbox{on}\ {\cal M}_\tau,
\end{equation*}
we have for all $\tau$
\begin{equation*}
\int_{\cG_d} |u|^2dx\le 2\tau^2 |\cG_d|+
2\int_{\cal M_\tau} (|u|-\tau)^2dx.
\end{equation*}
Let us take 
\begin{equation*}
\tau^2={1\over 4|\cG_d|} \int_{\cG_d} |u|^2 dx,
\end{equation*} 
i.e.
$\tau={1\over 2}\left(\overline{|u|^2}\right)^{1/2}$. Then for
this particular value of $\tau$ we obtain
\begin{equation}\label{E:(3)}
\int_{\cG_d} |u|^2 dx\le
4\int_{\cal M_\tau} (|u|-\tau)^2dx. 
\end{equation}
Assume first that
$\capa(\bcG_d\setminus {\cal M}_\tau)\ge\gamma \capa(\bcG_d)$.
Using \eqref{E:(3)} and applying Lemma \ref{L:Mazya1} to the function $(|u|-\tau)_+$, which
equals $|u|-\tau$ on ${\cal M}_\tau$ and $0$ on
$\cG_d\setminus {\cal M}_\tau$, we see that
\begin{equation*}
\capa(\bcG_d\setminus{\cal M}_\tau)\le 
\frac{C(\cG)\int_{{\cal M}_\tau}|\nabla(|u|-\tau)|^2 dx}
{|\cG_d|^{-1}\int_{\cG_d}|u|^2 dx}\le
\frac{C(\cG)\int_{\cG_d}|\nabla u|^2 dx}
{|\cG_d|^{-1}\int_{\cG_d}|u|^2 dx}\;,
\end{equation*}
where $C(\cG)$ is the same as in \eqref{E:cap-above}. Thus,
\begin{equation*}\label{E:(4-)}
\int_{\cG_d}|u|^2 dx \le\frac{C(\cG) |\cG_d|\int_{\cG_d}|\nabla u|^2 dx}
{\capa(\bcG_d\setminus{\cal M}_\tau)}
\le\frac{C(\cG) |\cG_d|\int_{\cG_d}|\nabla u|^2 dx}{\gamma\,\capa(\bcG_d)}
\end{equation*}
Note that  $|\cG_d|=|\cG| d^n$ and 
$\capa(\bcG_d)=\capa(\bcG) d^{n-2}$, where for $n=2$ 
the capacities of $\bcG=\bcG_1(0)$ and $\bcG_d=\bcG_d(y)$ are taken with respect 
to the discs $B_2(0)$ and $B_{2d}(y)$ respectively. Therefore we obtain
\begin{equation}\label{E:(4)}
\int_{\cG_d}|u|^2 dx 
\le\frac{C(\cG)|\cG| d^2}{\gamma\capa(\bcG)}\int_{\cG_d}|\nabla u|^2 dx.
\end{equation}

Now consider the opposite case
$\capa(\bcG_d\setminus {\cal M}_\tau)\le \gamma\, \capa(\bcG_d)$.
Then we can write
\begin{align*}
&\int_{\bcG_d} |u|^2 \cV(dx)\ge \int_{{\cal M}_\tau} |u|^2 \cV(dx)\ge
\tau^2 \cV({{\cal M}_\tau}) 
={1\over 4|\cG_d|}\int_{\cG_d}|u|^2dx\cdot \cV({\cal M}_\tau)\\
&\ge {1\over 4|\cG_d|}\int_{\cG_d}|u|^2dx\cdot \cV_\ga(\cG_d,\Om).
\end{align*}
Finally we obtain in this case
\begin{equation}\label{E:(5)}
\int_{\cG_d} |u|^2 dx\le 
\frac{4 |\cG_d|}{\cV_\ga(\cG_d,\Om))}\int_{\bcG_d} |u|^2 \cV(dx).
\end{equation}

The desired inequality \eqref{E:Mazya(2)} immediately follows from \eqref{E:(4)} and \eqref{E:(5)}
with $C_2(\cG)=\max\left\{C(\cG)|\cG|(\capa(\bcG))^{-1},4|\cG|\right\}$. $\square$

\ms
Now we will move to the proof of the sufficiency part in Theorem \ref{T:discr}. 

We will start with the following proposition which  gives 
a general (albeit complicated) sufficient condition for the discreteness of spectrum. 

\begin{proposition}\label{P:general-suff}
Given an operator $H_{\cV}$, let us assume that the following condition is satisfied:
there exists $\eta_0>0$ such that for every $\eta\in (0,\eta_0)$ we can find
$d=d(\eta)>0$ and $R=R(\eta)>0$, so that if $\cG_d$ satisfies  
$\bcG_d\cap (\Om\setminus B_R(0))\ne\emptyset$,  
then there exists $\ga=\ga(\cG_d, \eta)\in(0,1)$ such that
\begin{equation}\label{E:2-main-ineq}
\ga d^{-2}\ge \eta^{-1}\quad \text{and}\quad 
{d^{-n}}{\cV_\ga(\cG_d,\Om)}\ge  \eta^{-1}\;.
\end{equation}
Then the spectrum of $H_\cV$ is discrete.
\end{proposition}

{\bf Proof.} 
Recall that the discreteness of spectrum is equivalent to
the following condition: for every $\eta>0$ there exists $R>0$
such that \eqref{E:smalltails2} holds for every $u\in C_0^\infty(\Om)$.
This will be true if we establish that for every $\eta>0$ there exist
$R>0$ and
$d>0$ such that 
\begin{equation}\label{E:cube-tail}
\int_{\cG_d}|u|^2 dx\le \eta \left[\int_{\cG_d}|\nabla u|^2 dx +
\int_{\bcG_d} |u|^2\cV(dx)\right], 
\end{equation}
for all  $\cG_d$ such that $\bcG_d\cap(\Om\setminus B_R(0))\ne\emptyset$ 
and for all $u\in C^\infty(\bcG_d)$, such that  $u=0$ in a neighborhood of
$\bcG_d\setminus \Om$. 
Indeed, assume that 
\eqref{E:cube-tail} is true. Let us take a covering of  $\R^n$ by bodies $\bcG_d$ so that it has
a finite multiplicity $m=m(\cG)$ (i.e. at most $m$ bodies $\bcG_d$ can have non-empty intersection). 
Then, taking $u\in C_0^\infty(\Om)$ and summing up the estimates
\eqref{E:cube-tail} over all bodies $\cG_d$ with $\bcG_d\cap(\Om\setminus B_R(0))\ne\emptyset$, 
we obtain  \eqref{E:smalltails2} (hence \eqref{E:smalltails}) with $m\eta$ instead of $\eta$.

Now Lemma \ref{L:Mazya2} and the assumptions \eqref{E:2-main-ineq}
immediately imply \eqref{E:cube-tail} (with $\eta$ replaced by $C_2(\cG)\eta$). $\square$ 

\ms
Instead of requiring that the conditions of Proposition \ref{P:general-suff} are satisfied 
for all $\eta\in(0,\eta_0)$, it suffices to require it for a monotone sequence $\eta_k\to +0$.
We can also assume that $d(\eta_k)\to 0$ as $k\to+\infty$. Then, passing to a subsequence, 
we can assume that the sequence $\{d(\eta_k)\}$ is strictly decreasing.
Keeping this in mind, we can replace the dependence $d=d(\eta)$ by the inverse
dependence $\eta=g(d)$, so that $g(d)>0$ and $g(d)\to 0$ as $d\to +0$ (and here we can also
restrict ourselves to a sequence $d_k\to +0$). This leads to the following, 
essentially equivalent but more convenient reformulation of 
Proposition \ref{P:general-suff}:

\begin{proposition}\label{P:general-suff-d}
Given an operator $H_{\cV}$, assume that the following condition is satisfied:
there exists $d_0>0$ such that for every 
$d\in (0,d_0)$ we can find $R=R(d)>0$ and  $\ga=\ga(d)\in (0,1)$,
so that if 
$\bcG_d\cap (\Om\setminus B_R(0))\ne\emptyset$, then
\begin{equation}\label{E:2-main-ineq-d}
d^{-2}\ga\ge g(d)^{-1}\quad \text{and}\quad d^{-n}{\cV_\ga(\cG_d,\Om)}\ge g(d)^{-1},
\end{equation}
where  
$g(d)>0$ and $g(d)\to 0$ as $d\to +0$. Then the spectrum of $H_\cV$ is discrete.
\end{proposition}

{\bf Proof of Theorem \ref{T:discr}, part (ii).}
Instead of (ii) in Theorem \ref{T:discr} it suffices to prove  
the  (stronger) statement formulated in Remark \ref{R:other-suff}. 
So suppose that $\exists\; d_0>0$, $\exists\; c>0$, $\forall\; d\in(0,d_0)$,
$\exists\; R=R(d)>0$, $\exists\ga(d)\in (0,1)$, satisfying \eqref{E:gamma-cond1}, 
such that \eqref{E:inf-ineq}
holds for all $\cG_d$ with $\bcG_d\cap (\Om\setminus B_R(0))\ne\emptyset$.

Since the left hand side of \eqref{E:inf-ineq} is exactly
$d^{-n}\cV_{\ga(d)}(\cG_d,\Om)$, 
we see that
\eqref{E:inf-ineq} can be rewritten in the form 
\begin{equation*}
d^{-n}\cV_\ga(\cG_d,\Om)\ge cd^{-2}\ga(d),
\end{equation*}
hence we can apply Proposition \ref{P:general-suff-d} with $g(d)=c^{-1}d^2\ga(d)^{-1}$
to conclude that the spectrum of $H_\cV$ is discrete. $\square$

\section{A sufficiency precision example}
In this section we will prove Theorem \ref{T:precise}.
We will construct  a domain $\Om\subset \R^n$,  such that 
the condition \eqref{E:omega-cond} is satisfied with $\ga(d)=Cd^2$
(with an arbitrarily large $C>0$), and yet the spectrum of $-\De$ in $L^2(\Om)$
(with the Dirichlet boundary condition) is not discrete.  
This will show that the condition \eqref{E:gamma-cond1} is precise, so Theorem \ref{T:precise}
will be proved. 
We will assume for simplicity that $n\ge 3$.

\ms
We will use the following notations:

\begin{itemize}
\item{
$L^{(j)}$ is the spherical layer $\{x \in \R^{n}:\, \log j \leq |x|
\leq \log (j+1)\}$.
Its width is $\log(j+1)-\log j$ which is $<j^{-1}$ for all $j$ and 
equivalent to $j^{-1}$ for large $j$.
}
\item{
$\{Q_{k}^{(j)}\}_{k \geq 1}$ is a collection of closed cubes 
which form a tiling of $\R^n$ and have
edge length $\epsilon(n)\, j^{-1}$, where
$\epsilon(n)$ is a sufficiently small constant depending on $n$
(to be adjusted later).
}
\item{
$x_{k}^{(j)}$ is the center of $Q_{k}^{(j)}$.
}
\item{
$\{B_{k}^{(j)}\}_{k\geq 1}$ is the collection of closed balls centered
at $x_{k}^{(j)}$ with radii $\rho_{j}$ given by
\begin{equation*}
\om_n(n-2)\, \rho_{j}^{n-2}=C(\epsilon(n)/j)^{n},
\end{equation*}
where $\om_n$ is the area of the unit sphere $S^{n-1}\subset \R^n$ and $C$ is an
arbitrary constant. The last equality can be written as
\begin{equation}\label{UNO}
{\capa}(B_{k}^{(j)}) = C\, {\mes}\, Q_{k}^{(j)},
\end{equation}
where $\mes$ is the $n$-dimensional Lebesgue measure on $\R^n$.
Among the balls $B_{k}^{(j)}$ we will select a subcollection which consists of the balls 
with the additional property
$B_{k}^{(j)}\subset L^{(j)}$. We will  refer to these balls as \textit{selected} ones. 
We will denote selected balls by $\tB_{k}^{(j)}$. By an abuse of notation
we will not introduce special letter for the subscripts of the selected balls.
We will also denote by $\tQ^{(j)}_k$ the corresponding cubes $Q^{(j)}_k$,
so that 
\begin{equation*}
\tQ^{(j)}_k=Q^{(j)}_k\supset \tB_{k}^{(j)}.
\end{equation*}
}
\item{
$\Lambda^{(j)}=\bigcup_{k \geq 1}  \tB_{k}^{(j)}\subset L^{(j)}$.
\smallskip
}
\item{
$\Omega$ is the complement of $\cup _{j \geq 1} \Lambda^{(j)}$.
\smallskip
}
\item{
$B_{r}(P)$ is the closed ball with radius $r \leq 1$ centered at a
point $P$. We will make a more precise choice of $r$ later.
}
\end{itemize}

\begin{proposition} \label{P:1} The spectrum of  $-\Delta$ in
$\Omega$ (with the Dirichlet boundary condition) is not discrete.
\end{proposition}

{\bf Proof.} 
Let $j\ge 7$ and  $P
\in L^{(j)}$, i.e.
\begin{equation*}
\log j \leq |P| \leq \log(j+1).
\end{equation*}
Note that the ball $B_{r}(P)$ is a subset of the spherical layer
$\cup_{l\geq s \geq m} L^{(s)}$ if and only if
\begin{equation*}
\log m  \leq |P|-r \quad {\rm and}\quad |P|+r \leq \log(l+1).
\end{equation*}
Therefore, if
\begin{equation*}
\log m \leq \log j -r
\end{equation*}
and
\begin{equation*}
\log(j+1) +r \leq \log(l+1),
\end{equation*}
then $B_{r}(P) \subset \cup _{l \geq s \geq m} L^{(s)}$. The
last two inequalities can be written as
\begin{equation}\label{DUE}
m \leq j\, {e}^{-r}\quad {\rm and}\quad j+1 \leq (l+1)
{e}^{-r}.
\end{equation}
If we take, for example,
\begin{equation*}
m=[j/3]\quad {\rm and} \quad l=3j,
\end{equation*}
then, due to the inequality $j\ge 7$, we easily deduce that 
\begin{equation}\label{E:DUE.5}
B_{r}(P) \subset \bigcup_{[j/3]\leq s\leq 3j} L^{(s)}.
\end{equation}
Using \eqref{DUE}, the definition of $\Omega$ and subadditivity of  capacity, we obtain:
\begin{align*}
{\capa}(B_{r}(P) \setminus \Omega)&= 
{\capa}(B_{r}(P) \cap (\cup_{s\geq 1} \Lambda^{(s)}))\\
&\leq
\sum_{[j/3]\leq s \leq 3j}\ \sum_{k \geq 1} {\capa} (B_{r}(P)
\cap \tB_{k}^{(s)})\\
&\leq C \, \sum_{[j/3]\leq s \leq 3j}\ 
\sum_{\{k:B_{r}(P) \cap \tQ_{k}^{(s)} \neq \emptyset\}} 
{\mes}\, \tQ_{k}^{(s)}.
\end{align*}
It is easy to see that the multiplicity of the covering of $B_r(P)$ by the cubes
$\tQ_k^{(s)}$, participating in the last sum, is at most 2, provided $\epsilon(n)$
is chosen sufficiently small.
Hence,
\begin{equation}{\label{TRE}}
{\capa}(B_{r}(P) \setminus \Omega)\leq c(n)\, C \, r^{n}.
\end{equation}
On the other hand, we know that the discreteness of spectrum
guarantees that for every $r>0$
\begin{equation*}
\underset{|P| \to \infty}\liminf \> {\capa}(B_{r}(P) \setminus \Omega)
\geq \gamma(n) \, r^{n-2},
\end{equation*}
where $\gamma(n)$ is a constant depending only on $n$ (cf. Remark \ref{R:domains}). 
For sufficiently small $r>0$ this clearly contradicts \eqref{TRE}. 
$\square$

\begin{proposition}\label{P:2}
The domain  $\Omega$ satisfies
\begin{equation}\label{QUATTRO}
\underset{|P| \to \infty}{\liminf}{\capa}(B_{r}(P)
\setminus \Omega)\geq \delta(n)\, C\, r^{n},
\end{equation}
{\it where} $\delta(n)>0$ {\it depends only on} $n$.
\end{proposition}

{\bf Proof.} Let $\mu_{k}^{(s)}$ be the capacitary measure on
$\partial \tB_{k}^{(s)}$ 
(extended by zero to $\R^{n} \setminus \partial \tB_{k}^{(s)}$), 
and let $\epsilon_{1} (n)$ denote a sufficiently small
constant to be chosen later. We introduce the measure
\begin{equation*}
\mu= \epsilon_{1}(n) \sum_{k,s} \mu_{k}^{(s)},
\end{equation*}
where the summation here and below is taken over $k,s$ 
which correspond to the selected  balls $\tB_{k}^{(s)}$.
Taking $P\in L^{(j)}$, let us show that
\begin{equation}\label{CINQUE}
\int_{B_{r/2}(P)} \cE(x-y){d \mu(y)}\leq 1 \quad
{\rm on}\quad \R^{n},
\end{equation}
where $\cE(x)$ is given by \eqref{E:fund-sol}.
It suffices to verify (\ref{CINQUE}) for $x \in B_{r}(P)$,
because for $x \in \R^{n}\setminus B_{r}(P)$ this will follow
from the maximum principle.

Obviously, the potential in (\ref{CINQUE}) does not exceed
\begin{equation*}
\sum _{\{s,k:\tB_{k}^{(s)}\cap B_{r/2}(P)\ne\emptyset\}}
\epsilon_{1}(n)
\int_{\partial \tB_{k}^{(s)}} \cE(x-y)d \mu_{k}^{(s)}(y).
\end{equation*}
We divide this sum into two parts $\sum '$ and $\sum ''$, the first
sum being extended over all points $x_{k}^{(s)}$ with the
distance $\leq j^{-1}$ from $x$. Recalling that $x\in B_r(P)$ and using
\eqref{E:DUE.5}, we easily see that
the number of such points does not exceed
a certain constant $c_{1}(n)$. We define the constant
$\epsilon_{1}(n)$ by
\begin{equation*}
\epsilon_{1}(n)=(2 c_{1}(n))^{-1}.
\end{equation*}

Since $\mu_{k}^{(s)}$ is the capacitary measure, we have
\begin{equation*}
\sum {'}\ldots \leq \epsilon_{1}(n)\, c_{1}(n) = 1/2.
\end{equation*}
Furthermore, by (\ref{UNO})
\begin{align*}
\sum {''}\ldots &\leq c_{2}(n) \sum {''}
{{{\capa}(\tB_{k}^{(s)})}\over{|x-x_{k}^{(s)}|^{n-2}}}\
= c_{2}(n) \, C\ \sum {''}
{{{\mes}\,\tQ_{k}^{(s)}}\over{|x-x_{k}^{(s)}|^{n-2}}}\\
&\leq c_{3}(n) \ C\ \int_{B_{r}(P)}
{{dy}\over{|x-y|^{n-2}}} < c_{4}(n)\ C\ r^{2}.
\end{align*}
We can assume that
\begin{equation*}
r \leq (2 c_{4}(n) C)^{-1/2}
\end{equation*}
which implies $\sum '' \leq 1/2$. Therefore (\ref{CINQUE}) holds.

It follows that for large $|P|$ (i.e. for $P$ with $|P|\ge R=R(r)>0$), or, equivalently,
for large $j$,  we will have
\begin{align*}
{\capa}(B_{r}(P) \setminus \Omega) &\ge
 \sum _{\{s,k:\, \tB_{k}^{(s)} \subset B_{r/2}(P)\}}
 \epsilon_{1}(n)\mu_{k}^{(s)}(\partial \tB_{k}^{(s)})\\
&= \epsilon_{1}(n)  \sum _{{\{s,k:\, \tB_{k}^{(s)} \subset B_{r/2}(P)\}}} 
{\capa}(\tB_{k}^{(s)})\\
&= \epsilon_{1}(n)\ C\  \sum _{\{s,k:\, \tB_{k}^{(s)} \subset B_{r/2}(P)\}}
{\mes}\,Q_{k}^{(s)}\geq \delta(n) \ C\ r^{n}.
\end{align*}
This ends the proof of Proposition \ref{P:2}, hence of 
Theorem \ref{T:precise}. $\square$

\ms
\begin{remark}{\rm
Slightly modifying the construction given above, it is easy to provide an example of an
operator $H=-\De+V(x)$ with $V\in C^\infty(\R^n)$, $n\ge 3$, $V\ge 0$, such that the corresponding
measure $Vdx$ satisfies \eqref{E:inf-cond} with $\ga(d)=Cd^2$ and an arbitrarily large $C>0$,
but the spectrum of $H$ in $L^2(\R^n)$ is not discrete. 
So the condition \eqref{E:gamma-cond1} is precise
even in case of the Schr\"odinger operators with $C^\infty$ potentials.
}
\end{remark}

\section{Positivity of $H_\cV$}\label{S:positivity}

In this section we  prove Theorem \ref{T:positivity}. 

\ms
\textbf{Proof of Theorem \ref{T:positivity} (necessity).} 
Let us assume that the operator $H_{\cV}$ is strictly positive. This implies that the estimate
\eqref{E:mu-estimate} holds with some $\eta>0$ for every $\cG_d$ (with an arbitrary $d>0$) 
and every $u\in C_0^\infty(\cG_d\cap\Om)$.  But then we can use the arguments of Section \ref{S:nec} 
which lead to \eqref{E:main-nec-ineq-inf}, provided \eqref{E:eta-C2} is satisfied.  It will be satisfied if
$d$ is chosen sufficiently large.
$\square$

\ms
\textbf{Proof of Theorem \ref{T:positivity} (sufficiency).} 
Let us assume that there exist 
$d>0$, $\varkappa>0$ and $\ga\in(0,1)$ such that for every $\cG_d$
the estimate \eqref{E:pos-cond2} holds. Then 
by Lemma \ref{L:Mazya2},
for every $\cG_d$ and every $u\in C^\infty(\bcG_d)$, such that $u=0$ 
in a neighborhood of $\bcG_d\setminus \Om$, we have
\begin{equation*}
\int_{\cG_d}|u|^2 dx\le 
\frac{C_2(\cG)d^2}{\ga}\int_{\cG_d}|\nabla u|^2 dx +
\frac{C_2(\cG)d^n}{\varkappa}\int_{\bcG_d}|u|^2 \cV(dx).
\end{equation*}
Let us take a covering of $\R^n$ 
of finite multiplicity $N$ by bodies $\bcG_d$. 
It follows that for every $u\in C_0^\infty(\Om)$
\begin{equation*}
\int_{\Om}|u|^2 dx\le 
NC_2(\cG)d^2\max\left\{\frac{1}{\ga},\frac{d^{n-2}}{\varkappa}\right\}
\left(\int_{\Om}|\nabla u|^2 dx +
\int_{\Om}|u|^2 \cV(dx)\right),
\end{equation*}
which proves positivity of $H_{\cV}$. $\square$

\end{document}